\documentclass{amsart}
\usepackage{amssymb}
\usepackage{amscd, xypic}
\usepackage{times,amsmath,amssymb}

\newcommand\datver[1]{\def\datverp
{\par\boxed{\boxed{\text{Version: #1; Run: \today}}}}}
\datver{1.2; Revised: 04/28/2004}

\newcommand{\<}{\langle}
\renewcommand{\>}{\rangle}

\newcommand\maH{\mathcal H}

\newcommand\maF{\mathcal F}
\newcommand\maI{\mathcal I}

\newcommand\maC{\mathcal C}

\newcommand\maO{\mathcal O}

\newcommand\mfI{{\mathfrak I}}
\newcommand\kk{\mathfrak k}

\newcommand\per{\operatorname{per}}

\newcommand\CC{\mathbb C}
\newcommand\FF{\mathbb F}
\newcommand\GG{\mathbb G}
\newcommand\NN{\mathbb N}
\newcommand\QQ{\mathbb Q}

\newcommand\ZZ{\mathbb Z}

\newcommand\rg{\operatorname{reg}}
\newcommand\nil[1]{{\mathcal U}_{#1}}

\newcommand\pa{\partial}

\newcommand\CI{\mathcal{C}^\infty}
\newcommand\CIc{\mathcal{C}_c^\infty}
\newcommand\Tt{a_0 \otimes a_1 \otimes \ldots \otimes a_n}
\newcommand\ie{{\em i.e. }}

\newcommand\HH{\operatorname{HH}}
\newcommand\Hd{\operatorname{HH}}
\newcommand\Hc{\operatorname{HC}}

\newcommand\HP{\operatorname{HP}}
\newcommand\Hp{\operatorname{HP}}
\newcommand\Tr{\operatorname{Tr}}

\newcommand{\cohom}{\operatorname{H}}

\newcommand{\HeG}{\CIc(G)}
\newcommand{\Prim}{\operatorname{Prim}}
\newcommand{\an}{\operatorname{an}}
\newcommand{\mfk}{\mathfrak}
\newcommand\tHd{\operatorname{HH}^{\mathrm{top}}}

\newtheorem{theorem}{Theorem}[section]
\newtheorem{proposition}[theorem]{Proposition}
\newtheorem{corollary}[theorem]{Corollary}

\newtheorem{lemma}[theorem]{Lemma}

\theoremstyle{definition}
\newtheorem{definition}[theorem]{Definition}
\theoremstyle{remark}

\newtheorem{example}[theorem]{Example}

\author[V. Nistor]{Victor Nistor} \address{Pennsylvania State
University, University Park, PA 16802} \email{nistor@math.psu.edu}

\thanks{Nistor was partially supported by NSF Grant DMS-0200808.
Manuscripts available from {\bf http:{\scriptsize//}www.math.psu.edu}}

\begin{document}

\dedicatory\datverp

\title[Homology of Hecke algebras]{A non-commutative geometry approach
to the representation theory of reductive $p$-adic groups:
Homology of Hecke algebras, a survey and some new results}


\begin{abstract}\
We survey some of the known results on the relation between the
homology of the {\em full} Hecke algebra of a reductive $p$-adic
group $G$, and the representation theory of $G$. Let us denote by
$\CIc(G)$ the full Hecke algebra of $G$ and by $\Hp_*(\CIc(G))$
its periodic cyclic homology groups. Let $\hat G$ denote the
admissible dual of $G$. One of the main points of this paper is
that the groups $\Hp_*(\CIc(G))$ are, on the one hand, directly
related to the topology of $\hat G$ and, on the other hand, the
groups $\Hp_*(\CIc(G))$ are explicitly computable in terms of $G$
(essentially, in terms of the conjugacy classes of $G$ and the
cohomology of their stabilizers). The relation between
$\Hp_*(\CIc(G))$ and the topology of $\hat G$ is established as
part of a more general principle relating $\Hp_*(A)$ to the
topology of $\Prim(A)$, the primitive ideal spectrum of $A$, for
any finite typee algebra $A$.  We provide several new examples
illustrating in detail this principle. We also prove in this paper
a few new results, mostly in order to better explain and tie
together the results that are presented here. For example, we
compute the Hochschild homology of $\maO(X) \rtimes \Gamma$, the
crossed product of the ring of regular functions on a smooth,
complex algebraic variety $X$ by a finite group $\Gamma$. We also
outline a very tentative program to use these results to construct
and classify the cuspidal representations of $G$. At the end of
the paper, we also recall the definitions of Hochschild and cyclic
homology.
\end{abstract}

\maketitle
\tableofcontents

\section*{Introduction\label{Sec.I}}

To motivate the results surveyed in this paper, let us look at the
following simple example. Precise definitions will be given below.
Let $G$ be a finite group and $A := \CC[G]$ be its complex group
algebra. Then $A$ is a finite dimensional, semi-simple complex
algebra, and hence $A \simeq \oplus_{j=1}^d M_{n_j}(\CC)$. (This
is an elementary result that can be found in \cite{Lang}; see also
\cite{Serre}). The Hochschild homology of $A$ is then, on the one
hand,
\begin{equation}
    \HH^0(A) \simeq \oplus_{j=1}^d \HH^0(M_{n_j}(\CC))
    \simeq \CC^d.
\end{equation}
On the other hand, $\HH^0(A)$ is the space of traces on $A$, and
hence it identifies with the space of class functions on $G$. Let
$\<G\>$ denote the set of conjugacy classes of $G$ and $\#S$
denote the number of elements in a set $S$. Then $\HH^0(A)$ has
dimension $\#\<G\>$. In other words,

\begin{proposition}[Classical]\label{prop.class}
Let $G$ be a finite group. Then $G$ has as many (equivalence classes
of) irreducible, complex representations as conjugacy classes.
\end{proposition}

One of our goals was to investigate to what extent Proposition
\ref{prop.class} extends to other groups. It is clear that the
formulation of any possible extension of Proposition
\ref{prop.class} will depend on the class of groups considered and
will not be as simple as in the finite group case. Moreover, this
question will not be answered in a few papers and is more of a
program (going back to Gelfand, Langlands, Manin, and other
people) than an explicit question. Nevertheless, something from
Proposition \ref{prop.class} does remain true in certain cases. An
example is the theory of characters for compact Lie groups.

In this paper, we will investigate a possible analog of
Proposition \ref{prop.class} for the case of a reductive $p$--adic
group $G$. Recall that a $p$--adic group $G = \GG(\FF)$ is the set
of $\FF$--rational points of a linear algebraic group $\GG$
defined over a non-archemedean, non-discrete, locally compact
field $\FF$ of characteristic zero (so $\FF$ is a finite algebraic
extension of the field $\QQ_l$ of $l$-adic numbers, for some prime
$l$). A vague formulation of our main result is as follows.

Let $\HP_k(A)$ denote the periodic cyclic homology groups of the
algebra $A$. Also, let $\HeG$ denote the space of compactly
supported, locally constant functions on $G$ and let $\hat G$
denote the admissible dual of $G$ with the Jacobson topology. Then
the results of \cite{BaumNistor0, BaumNistor, KNS, Nistor29} give
the following result that will be made more precise below.

\begin{theorem}\label{thm.vague}\
The groups $\HP_j(\HeG)$ are explicitly determined by the geometry
of the conjugacy classes of $G$ and the cohomology of their
stabilizers and they are (essentially) isomorphic to the singular
cohomology groups of $\hat G$.
\end{theorem}

One of the main purposes of this paper is to explain the above
theorem. This theorem is useful especially because it is much
easier to determine the groups $\HP_j(\HeG)$ (and hence, to a
large extent, the algebraic cohomology of $\hat G$) than it is to
determine $\hat G$ itself. Moreover, we will briefly sketch a plan
to say more about the actual structure of $\hat G$ using the
knowledge of the topology of $\hat G$ acquired from the
determination of $\HP_j(\HeG)$ in terms of the geometry of the
conjugacy classes of $G$ and the cohomology of their stabilizers.
See also \cite{BHP} for a survey of the applications of
non-commutative geometry to the representation theory of reductive
$p$-adic groups.

The paper is divided into two parts. The first part, consisting of
Sections \ref{sec1}--\ref{sec4} is more advanced, whereas the last
three sections review some basic material. In Section \ref{sec1}
we review the basic result relating the cohomology of the maximal
spectrum of a commutative algebra $A$ to its periodic cyclic
homology groups $\Hp_*(A)$. The relation between forms on
$\operatorname{Max}(A)$ and the Hochschild homology groups
$\Hd_*(A)$ are also discussed here. These are basic results due to
Connes \cite{ConnesNCG}, Feigin and Tsygan \cite{FeiginTsygan},
and Loday and Quillen \cite{LodayQuillen}. We also discuss the
Excisition principle in periodic cyclic homology
\cite{CuntzQuillen} and it's relation with $K$-theory. In Section
\ref{sec2} we discuss generalizations of these results to finite
type algebras, a class of algebras directly relevant to the
representation theory of $p$-adic groups that was introduced in
\cite{KNS}. We also use these results to compute the periodic
cyclic homology and the Hochschild homology of several typical
examples of finite type algebras. In the following section,
Section \ref{sec3}, we introduce spectrum preserving morphisms,
which were shown in \cite{BaumNistor} to induce isomorphisms on
periodic cyclic homology. This then led to a determination of the
periodic cyclic homology of Iwahori-Hecke algebras in that paper.
In Section \ref{sec4}, we recall the explicit calculation of the
Hochschild and periodic cyclic homology groups of the full Hecke
algebra $\CIc(G)$. The last three sections briefly review for the
benefit of the reader the definitions of Hochschild homology,
cyclic and periodic cyclic homology, and, respectively, the Chern
character.

This first part of the paper follows fairly closely the structure
of my talk given at the conference ``Non-commutative geometry and
number theory'' organized by Yuri Manin and Matilde Marcolli, whom
I thank for their great work and for the opportunity to present my
results. I have included, however, some new results, mostly to
better explain and illustrate the results surveyed.

\section{Periodic cyclic homology versus singular
cohomology\label{sec1}}

Let us discuss first to what extent the periodic cyclic homology
groups $\Hp(\HeG)$ are related to the topology -- more precisely
to the singular cohomology -- of $\hat G$. In the next section, we
will discuss this again in the more general framework of ``finite
type algebras'' (Definition \ref{def.FT}). The definitions of the
homology groups considered in this paper and of the Connes-Karoubi
character are recalled in the last three sections of this paper.

One of the main goals of non-commutative geometry is to generalize
the correspondence (more precisely, contravariant equivalence of
categories)
\begin{equation}\label{eq.NCG}
    \text{ \em ``Space'' } X \;\; \leftrightarrow \maF(X) := \text{ \em
    ``the algebra of functions on $X$'' }
\end{equation}
to allow for non-commutative algebras on the right--hand side of
this correspondence. This philosophy was developed in many papers,
including \cite{BBN, BPChern, BenameurNistor1, BrylinskiNistor,
ConnesNCG, ConnesBook, Karoubi, Kontsevich, Loday, LodayQuillen,
Manin1, Manin2, ManinMarcolli, Nistor5, Nistor8}, to mention only
some of the more recent ones. The study of the $K$-theory of
$C^*$-algebras, a field on its own, certainly fits into this
philosophy. The extension of the correspondence in Equation
\ref{eq.NCG} would lead to methods to study (possibly)
non-commutative algebras using our geometric intuition. In
Algebraic geometry, this philosophy is illustrated by the
correspondence (\ie contravariant equivalence of categories)
between affine algebraic varieties over a field $\kk$ and
commutative, reduced, finitely generated algebras over $\kk$. In
Functional analysis, this principle is illustrated by the
Gelfand--Naimark equivalence between the category of compact
topological spaces and the category of commutative, unital
$C^*$-algebras. In all these cases, the study of the ``space''
then proceeds through the study of the ``algebra of functions on
that space.'' For this approach to be useful, one should be able
to define many invariants of $X$ in terms of $\maF(X)$ alone,
preferably without using the commutativity of $\maF(X)$.

It is a remarkable fact that one can give completely algebraic
definitions for $\Omega^q(X)$, the space of differential forms on $X$
(for suitable $X$) just in terms of $\maF(X)$. Even more remarkable is
that the singular cohomology of $X$ (again for suitable $X$) can be
defined in purely algebraic terms using only $\maF(X)$. In these
definitions, we can then replace $\maF(X)$ with a non-commutative
algebra $A$. Let us now recall these results.

We denote by $\HH_j(A)$ the Hochschild homology groups of an
algebra $A$ (see Section \ref{sec5} for the definition). As we
will see below, for applications to representation theory we are
mostly interested in the algebraic case, so we state those first
and then we state the results on smooth manifolds. We begin with a
result of Loday-Quillen \cite{LodayQuillen}, in this result
$\maF(X) = \maO(X)$, the ring of regular (\ie polynomial)
functions on the algebraic variety $X$.

\begin{theorem}[Loday-Quillen]\label{thm.LQ1}\
Let $X$ be a smooth, complex, affine algebraic variety. Then
\begin{equation*}
    \HH_j(\maO(X)) \simeq \Omega^j(X),
\end{equation*}
the space of algebraic forms on $X$.
\end{theorem}

A similar results holds when $X$ is a smooth compact manifold and
$\maF(X) = \CI(X)$ is the algebra of smooth functions on $X$
\cite{ConnesNCG}. See also the Hochschild--Kostant--Rosenberg
paper \cite{HKR}.

\begin{theorem}[Connes]\label{thm.Connes1}\
Let $X$ be a compact, smooth manifold. Then
\begin{equation*}
    \HH_j(\CI(X)) \simeq \Omega^j(X),
\end{equation*}
the space of smooth forms on $X$.
\end{theorem}

These results extend to recover the singular cohomology of
(suitable) spaces, as seen in the following two results due to
Feigin-Tsygan \cite{FeiginTsygan} and Connes \cite{ConnesNCG}. For
any functor $F_j$ [respectively, $F^j$], we shall denote by
$F_{[j]} = \oplus_{k}F_{j + 2k}$ [respectively, $F^{[j]} =
\oplus_{k}F^{j + 2k}$]. This will mostly be used for $F^j(X) =
\cohom^j(X)$, the singular cohomology of $X$.

\begin{theorem}[Feigin-Tsygan]\label{thm.FT1}\
Let $X$ be a complex, affine algebraic variety and $\maO(X)$ be
the ring of regular (\ie polynomial) functions on $X$. Then
\begin{equation*}
    \Hp_j(\maO(X)) \simeq \cohom^{[j]}(X).
\end{equation*}
\end{theorem}

For smooth algebraic varieties, this result follows from the
Loday-Quillen result on Hochschild homology mentioned above. See
\cite{KNS} for a proof of this theorem that proceed by reducing it
to the case of smooth varieties. For smooth manifolds, the result
again follows from the corresponding result on Hoschschild
homology.

\begin{theorem}[Connes]\label{thm.Connes2}\
Let $X$ be a compact, smooth manifold and $\CI(X)$ be the algebra
of smooth functions on $X$. Then
\begin{equation*}
    \Hp_j(\CI(X)) \simeq \cohom^{[j]}(X).
\end{equation*}
\end{theorem}

These results are already enough justification for declaring
periodic cyclic homology to be the ``right'' extension of singular
cohomology for the category (suitable) spaces to suitable
categories of algebras. However, the most remarkable result
justifying this is the ``Excision property'' in periodic cyclic
homology, a breakthrough result of Cuntz and Quillen
\cite{CuntzQuillen}.

\begin{theorem}[Cuntz-Quillen]\label{Theorem.Excision}\
Any two-sided ideal $J$ of an algebra $A$ over a characteristic
$0$ field gives rise to a periodic six-term exact sequence

\begin{equation}
\xymatrix{
\Hp_0(J)            \ar[r] & \Hp_0(A)  \ar[r] &
\Hp_0(A/J)          \ar[d]^{\pa} \\ \Hp_1(A/J) \ar[u]^{\pa} &
\Hp_1(A)                \ar[l] & \Hp_1(J). \ar[l] }
\end{equation}
\end{theorem}

A similar results holds for Hochschild and cyclic homology,
provided that the ideal $J$ is an $H$-unital algebra in the sense
of Wodzicki \cite{Wodzicki}, see Section \ref{sec5}. We shall refer
to the following result of Wodzicki as the ``Excision principle in
Hochschild homology.''

\begin{theorem}[Wodzicki]\label{thm.Wodz}\
Let $J \subset A$ be a $H$-unital ideal of a complex algebra $A$.
Then there exists a long exact sequence
\begin{multline*}
    0 \leftarrow \HH_0(A/J) \leftarrow HH_0(A) \leftarrow \HH_0(J)
    \stackrel{\pa}{\leftarrow} \HH_1(A/J) \\ \leftarrow \HH_k(A/J)
    \leftarrow HH_k(A) \leftarrow \HH_k(J)
    \stackrel{\pa}{\leftarrow} \HH_{k+1}(A/J) \leftarrow \ldots
\end{multline*}
The same result remains valid if we replace Hochschild homology with
cyclic homology.
\end{theorem}

Also, there exist excision results for topological algebras
\cite{BL, Cuntz:m}. An important part of the proof of the Excision
property is to provide a different definition of cyclic homology in
terms of $X$-complexes. Then the proof is ingeniously reduced to
Wodzicki's result on the excision in Hochschild homology, using also
an important theorem of Goodwillie \cite{Goodwillie} that we now
recall.

\begin{theorem}[Goodwillie]\label{theorem.Goodwillie}\
If $I \subset A$ is a nilpotent two-sided ideal, then the quotient
morphism $A \to A/I$ induces an isomorphism $\Hp_*(A) \to
\Hp_*(A/I)$. In particular, $\Hp_*(I) = 0$ whenever $I$ is nilpotent.
\end{theorem}

One of the main original motivations for the study of cyclic
homology was the need for a generalization of the classical Chern
character $Ch : K^j(X) \to \cohom^{[j]}(X)$ \cite{ConnesNCG,
ConnesBook, Karoubi, Loday}. Indeed, an extension is obtained in
the form of the Connes-Karoubi character
\begin{equation}
    Ch : K_j(A) \to \HP_j(A).
\end{equation}
It is interesting then to notice that excision in periodic cyclic
homology is compatible with excision in $K$-theory, which is seen from
the following result \cite{Nistor15}, originally motivated by
questions in the analysis of elliptic operators (more precisely, Index
theory).

\begin{theorem}[Nistor]\
Let $I \subset A$ the a two-sided ideal of a complex algebra $A$.
Then the following diagram commutes
\begin{equation*}
\xymatrix{
K_1(I)        \ar[d] \kern-1pt\ar[r]\kern-1pt & K_1(A)
\ar[d] \kern-1pt\ar[r]\kern-1pt & K_1(A/I)            \ar[d]
\kern-1pt\ar[r]^{\pa}\kern-1pt& K_0(I)          \ar[d]
\kern-1pt\ar[r]\kern-1pt     & K_0(A)          \ar[d]
\kern-1pt\ar[r]\kern-1pt     & K_0(A/I)            \ar[d]      \\
\Hp_1(I)            \kern-1pt\ar[r]\kern-1pt        & \Hp_1(A)
\kern-1pt\ar[r]\kern-1pt        & \Hp_1(A/I)
\kern-1pt\ar[r]^{\pa}\kern-1pt  & \Hp_0(I)
\kern-1pt\ar[r]\kern-1pt        & \Hp_0(A)
\kern-1pt\ar[r]\kern-1pt        & \Hp_0(A/I).                 }
\end{equation*}
\end{theorem}

Let $X$ be a complex, affine algebraic variety, $\maO(X)$ the ring
of polynomial functions on $X$, $X^{\an}$ the underlying locally
compact topological space, and $Y \subset X$ be a subvariety. Let
$I \subset \maO(X)$ be the ideal of functions vanishing on $Y$.
Then the above theorem shows, in particular, that the periodic six
term exact sequence of periodic cyclic homology groups associated
to the exact sequence
\begin{equation}
    0 \longrightarrow I \longrightarrow \maO(X) \longrightarrow
    \maO(Y) \longrightarrow 0,
\end{equation}
of algebras by the Excision principle is obtained from the long
exact sequence in singular cohomology of the pair $(X^{\an},
Y^{\an})$ by making the groups periodic of period two \cite{KNS,
Nistor15}. The same result holds true for the exact sequence of
algebras associated to a closed submanifold $Y$ of a smooth
manifold $X$.

\section{Periodic cyclic homology and $\hat G$\label{sec2}}

Let $A$ be an arbitrary complex algebra. The kernel of an
irreducible representation of $A$ is called a {\em primitive}
ideal of $A$. We shall denote by $\Prim(A)$ the {\em primitive
ideal spectrum} of $A$, consisting of all primitive ideals of $A$.
We endow $\Prim(A)$ with the Jacobson topology. Thus, a set $V
\subset \Prim(A)$ is open if, and only if, $V$ is the set of
primitive ideals not containing some fixed ideal $I$ of $A$. We have
\begin{equation*}
    \Prim(\maO(X)) =: \operatorname{Max}(\maO(X)) = X,
\end{equation*}
the set of maximal ideals of $A = \maO(X)$ with the Zariski
topology. If $A = \CI(X)$, where $X$ is a smooth compact manifold,
then again $\Prim(A) = \operatorname{Max}(A) = X$ with the usual (\ie
locally compact, Hausdorff) topology on $X$.
We are interested in the primitive ideal spectra of algebras because
\begin{equation}\label{eq.Bernstein}
    \Prim(\CIc(G)) = \hat G,
\end{equation}
a deep result due to Bernstein \cite{Deligne}. For the purpose of this
paper, we could as well take $\Prim(\CIc(G))$ to be the actual
definition of $\hat G$.

In view of the results presented in the previous section, it is
reasonable to assume that the determination of the groups
$\HP_j(\CIc(G))$ would give us some insight into the topology of $\hat
G$. I do not know any general result relating the singular cohomology
of $\hat G$ to $\HP_j(\CIc(G))$, although it is likely that they
coincide. Anyway, due to the fact that the topology on $\hat G$ is
highly non-Hausdorff topology, it is not clear that the knowledge of
the groups $\cohom^j(\hat G)$ would be more helpful than the knowledge
of the groups $\HP_j(\CIc(G))$.

For reasons that we will explain below, it will be convenient in what
follows to work in the framework of ``finite typee algebras''
\cite{KNS}. All our rings have a unit (\ie they are unital), but
the algebras are not required to have a unit.

\begin{definition}\label{def.FT}\
Let $\kk$ be a finitely generated commutative, complex ring. A
finite typee $\kk$-algebra is a $\kk$-algebra that is a finitely
generated $\kk$-module.
\end{definition}

The study of $\hat G$ as well as that of $A = \HeG$ reduces to the
study of finite typee algebras by considering the connected components
of $\hat G$ and their commuting algebras, in view of some results of
Bernstein \cite{Deligne} that we now recall. Let $D \subset \hat G$ be
a connected component of $\hat G$. Then $D$ corresponds to a cuspidal
representation $\sigma$ of a Levi subgroup $M \subset G$.  Let $M_0$
be the subgroup of $M$ generated by the compact subgroups of $M$ and
$H_D$ be the representation of $G$ induced from the restriction of
$\sigma$ to $M_0$. The space $H_D$ can be thought of as the
holomorphic family of induced representations of
$\operatorname{ind}_{M}^G(\sigma \chi)$, where $\chi$ ranges through
the caracters of $M/M_0$. Let $A_D$ be the algebra of
$G$-endomorphisms of $H_D$. This is Bernstein's celebrated ``commuting
algebra.'' The annihilator of $H_D$ turns out to be a direct summand
of $\CIc(G)$ with complement the two-sided ideal $\CIc(G)_D \subset
\CIc(G)$. Then the category of modules over $\CIc(G)_D$ is equivalent
to the category of modules over $A_D$. Our main reason for introducing
finite typee algebras is that the algebra $A_D$ is a unital finite typee
algebra and
\begin{equation}\label{eq.AD}
    D = \Prim(A_D).
\end{equation}
Moreover,
\begin{equation}\label{eq.DirectSum}
    \CIc(G)) = \oplus_{D} \CIc(G)_D.
\end{equation}
To get consequences for the periodic cyclic homology, we shall
need the following result.

\begin{proposition}\label{prop.DirectSum}\
Let $B$ be a complex algebra such that $B \simeq \oplus_{n \in
\NN} B_n$ for some two-sided ideals $B_n \subset B$. Assume that
$B = \cup e_kBe_k$, for a sequence of idempotents $e_k \in B$ and
that $\HH_q(B) = 0$ for $q > N$, for some given $N$. Then
\begin{equation*}
    \Hp_q(B) \simeq \oplus_{n \in \NN} \Hp_q(B_n).
\end{equation*}
\end{proposition}

\begin{proof}\
The algebra $B$ and each of the two-sided ideals $B_n$ are
$H$-unital. Then
\begin{equation*}
    \HH_q(B) \simeq \oplus_{n \in \NN} \HH_q(B_n),
\end{equation*}
by Wodzicki's excision theorem and the continuity of Hochschild
homology (\ie the compatibility of Hochschild homology with inductive
limits). Then $\HH_q(B_n) = 0$ for $q > N$ and any $n$.  Therefore
$\Hp_{k}(B_n) = \Hc_k(B_n)$ for $k \ge N$, by the SBI-long exact
sequence (this exact sequence is recalled in Section
\ref{sec6}). Unlike periodic cyclic homology, cyclic homology is
continuous (\ie it is compatible with inductive limits). Using again
Wodzicki's excision theorem, we obtain
\begin{equation*}
    \Hp_{k} = \Hc_k(B) \simeq \oplus_{n \in \NN} \Hc_k(B_n) =
    \oplus_{n \in \NN} \Hp_{k}(B_n),
\end{equation*}
for $k \ge N$. This completes the proof.
\end{proof}

The above discussion gives the following result mentioned in
\cite{KNS} without proof. For the proof, we shal also use Theorem
\ref{Theorem.Gen}, which implies, in particular, that $\HH_q(\CIc(G))$
vanishes for $q$ greater than the split rank of $G$. This is, in fact,
a quite non-trivial property of $\CIc(G)$ and of the finite type
algebras $A_D$, as we shall see below in Example \ref{ex.nb}.

\begin{theorem}\label{thm.DS}\
Let $D$ be the set of connected components of $\hat G$, then
\begin{equation*}
    \Hp_q(\CIc(G)) = \oplus_{D} \Hp_q(\CIc(G)_D) \simeq
    \oplus_{D} \Hp_q(A_D).
\end{equation*}
\end{theorem}

\begin{proof}\ The first part follows directly from the results
above, namely from Equation \eqref{eq.DirectSum}, Proposition
\ref{prop.DirectSum}, and Theorem \ref{Theorem.Gen} (which implies
that $\HH_q(\CIc(G)) = 0$ for $q > N$, for $N$ large).

To complete the proof, we need to check that
\begin{equation}
    \Hp_q(\CIc(G)_D) \simeq \Hp_q(A_D).
\end{equation}
Let $e_k$ be a sequence of idempotents of $\CIc(G)$ corresponding
to a basis of neighborhood of the identity of $G$ consisting of
compact open subgroups. Then, for $k$ large, the unital algebra
$e_k \CIc(G) e_k$ is Morita equivalent to $A_D$, an imprimitivity
module being given by $e_k H_D$. (Recall from above that $H_D$ is
the induced representation from a cuspidal representation of
$M_0$, where $M$ is a Levi subfactor defining the connected
component $D$ and $M_0$ is the subgroup of $M$ generated by its
compact subgroups.)

In particular, $\HH_q(e_k \CIc(G)_d e_k)$ vanish for $q$ large.
The same argument as above (using Theorem \ref{Theorem.Gen} below)
shows that $\Hp_q(\CIc(G)_D) \simeq \Hp_q(e_k \CIc(G)_D e_k)$, for
$k$ large. The isomorphism $\Hp_q(\CIc(G)_D) \simeq \Hp_q(A_D)$
then follows from the invariance of Hochschild homology with
respect to Morita equivalence.
\end{proof}

For suitable $G$ and $D$,
\begin{equation}
    A_D = H_q,
\end{equation}
that is, the commuting algebra of $D$ is the Iwahori-Hecke algebra
associated to $G$ (or to its extended affine Weyl group),
\cite{Borel}. The periodic cyclic homology groups of $H_q$ were
determined in \cite{BaumNistor0, BaumNistor}, and will be recalled in
Section \ref{sec.BC}.

In view of Equations \eqref{eq.AD} and \eqref{eq.DirectSum} and of
the Theorem \ref{thm.DS}, we see that in order to relate the
groups $\Hp_*(\CIc(G))$ to the topology of $\hat G$, it is enough
to relate $\Hp_*(A)$ to the topology of $\Prim(A)$ for an
arbitrary finite typee algebra. For the rest of this and the
following section, we shall therefore concentrate on finite typee
algebras and their periodic cyclic homology.

If $\mfI \subset A$ is a primitive ideal of the finite typee
$\kk$-algebra $A$, then the intersection $\mfI \cap Z(A)$ is a
maximal ideal of $Z(A)$. The resulting map
\begin{equation*}
\Theta : \Prim(A) \to \operatorname{Max}(Z(A))
\end{equation*}
is called the {\em central character} map. We similarly obtain a
map $\Theta : \Prim(A) \to \operatorname{Max}(\kk)$, also called
the central character map.

The topology on $\Prim(A)$ and the groups $\HP_j(A)$ are related
through a spectral sequence whose $E^2$ are given by the singular
cohomology of various strata of $\Prim(A)$ that are better behaved
than $\Prim(A)$ itself, Theorem \ref{Theorem.spectral.sq} below.  To
state our next result, due to Kazhdan, Nistor, and Schneider
\cite{KNS}, we need to introduce some notation and definitions.

We shall use the customary notation to denote by $0$ the ideal
$\{0\}$.  Recall that an algebra $B$ is called {\em semiprimitive}
if the intersection of all its primitive ideals is $0$. Also, we
shall denote by $Z(B)$ the center of an algebra $B$. We shall need
the following definition from \cite{KNS}

\begin{definition}\label{Def.abelian} \
A finite decreasing sequence
\begin{equation*}
    A = \mfk I_0 \supset \mfk I_1 \supset \ldots \supset
    \mfk I_{n-1} \supset \mfk I_n
\end{equation*}
of two-sided ideals of a unital finite type algebra $A$ is an {\em
abelian filtration} if and only if the following three conditions
are satisfied for each $k$:

(i) The quotient $A/\mfk I_k$ is semiprimitive.

(ii) For each maximal ideal $\mfk p \subset Z_k := Z(A/\mfk I_k)$ not
containing $I_{k-1} := Z_k \cap (\mfk I_{k-1}/\mfk I_k)$, the
localization $(A/{\mfk I}_k)_{\mfk p}=(Z_k \setminus \mfk
p)^{-1}(A/{\mfk I}_k)$ is an Azumaya algebra over $(Z_k)_{\mfk p}$,

(iii) The quotient $(\mfk I_{k-1}/{\mfk I_k})/I_{k-1}(A/\mfk I_k)$
is nilpotent and $\mfk I_n$ is the intersection of all primitive
ideals of $A$.
\end{definition}

Consider an abelian filtration $(\mfk I_k)$, $k = 0, \ldots, n$ of a
finite type $\kk$-algebra $A$. Then, for each $k$, the center $Z_k$ of
$A/{\mfk I}_k$ is a finitely generated complex algebra, and hence it
is isomorphic to the ring of regular functions on an affine, complex
algebraic variety $X_k$. Let $Y_k \subset X_k$ be the subvariety
defined by $I_{k-1} := Z_k \cap ({\mfk I}_{k-1}/{\mfk I}_{k})$. For
any complex algebraic variety $X$, we shall denote by $X^{\an}$ the
topological space obtained by endowing the set $X$ with the locally
compact topology induced by some embedding $X \subset \CC^N$, $N$
large. We shall call $X^{\an}$ the analytic space underlying $X$. For
instance, we shall refer to $X_k^{\an}$ and $Y_k^{\an}$ as the {\em
analytic spaces} associated to the filtration $({\mathfrak I}_k)$ of
$A$.  Then we have the following result from \cite{KNS}

\begin{theorem}[Kazhdan-Nistor-Schneider]\
\label{Theorem.spectral.sq} If $Y_p^{\an} \subset X_p^{\an}$ are
the analytic spaces associated to an abelian filtration of a
finite type algebra $A$, then there exists a natural spectral
sequence with
\begin{equation*}
    E^1_{-p,q} = \cohom^{[q-p]}(X_p^{\an},Y_p^{\an})
\end{equation*}
convergent to $\Hp_{q-p}(A)$.
\end{theorem}

If the algebra $A$ in the above theorem is commutative, then any
decreasing filtration of $A$ (by radical ideals) is abelian, which
explains the terminology ``abelian filtration.'' Moroever, in this
case our spectral sequence reduces to the spectral sequence in
singular cohomology associated to the filtration of a space by
closed subsets.

Every finite type algebra has an abelian filtration. Indeed, we
can take ${\mathfrak I}_k$ to be the intersection of the kernels
of all irreducible representations of dimension at most $k$. This
filtration is called the {\em standard filtration}.

\begin{example}\label{ex1}\
Let $a_1, a_2, \ldots, a_l \in \CC$ be distinct points and $v_1,
v_2, \ldots, v_l \in \CC^2$ be non-zero (column) vectors.  Let
$\CC[x] = \maO(\CC)$ denote the algebra of polynomials in one
variable and define
\begin{equation}\label{eq.def.A}\
    A_1 := \{P \in M_2(\CC[x]),\, F(a_j)
    v_j \in \CC v_j \}.
\end{equation}
Then $A_1$ is a finite type algebra with center $Z=Z(A_1)$, the
subalgebra of matrices of the form $P E_2$, where $P \in \CC[x]$
and $E_2$ is the identity matrix in $M_2(\CC)$.

We filter $A_1$ by the ideals
\begin{equation*}
    \mfI_1 = \{P \in A_1,\, F(a_j) = 0 \}
\end{equation*}
and $\mfI_2 = 0$. We have
\begin{equation*}
    A_1/\mfI_1\simeq \CC^{2k}
\end{equation*}
a semi-primitive (\ie reduced) commutative algebra with center
\begin{equation*}
    Z_1 := Z(A_1/\mfI_1) = A_1/\mfI_1 = \CC^{2k}.
\end{equation*}
Then $I_0 := Z(A_1/\mfI_1)
\cap \mfI_0 = A_1$ and hence no maximal ideal $\mfk p$ contains
$I_0$. On the other hand, the algebra $A_1/\mfI_1$ is an Azumaya
algebra, so Condition (ii) of Definition \ref{Def.abelian} is
automatically satisfied for $k = 1$.  Similarly, the quotient
$(\mfI_{0}/{\mfI_1})/I_{0}(A_1/\mfI_1) = 0$, and hence it is
nilpotent. Hence Condition (iii) of Definition \ref{Def.abelian} is
also satisfied for $k = 1$. The space $X_1$ consists of $2l$ points
(two copies of each $a_j$) and $Y_1$ is empty.

Next, $Z_2 = Z \simeq \CC[x]$ and $A_1 = A_1/\mfI_2$ is
semi-primitive. The ideal $I_1 := Z_2 \cap \mfI_1$ consists of the
polynomials vanishing at the given values $a_j$. In this case,
$X_2 = \CC$ and $Y_2 = \{a_1, \ldots, a_l\}$. The Conditions
(i--iii) of Definition \ref{Def.abelian} are easily checked for $k
= 2$. This also follows from the fact that $\maI_j$ is the
standard filtration of $A_1$.

The spectral sequence of Theorem \ref{Theorem.spectral.sq} then
becomes $E^{1}_{-p, q}=0$, unless $p = 1$ or $p = 2$, in which case we
get
\begin{equation*}
    E^{1}_{-1, q} = \cohom^{[q-1]}(X_1, Y_1) =
\begin{cases}
\CC^{2l} & \text{if } q \text{ is odd}\\ 0 & \text{ otherwise}
\end{cases}
\end{equation*}
and
\begin{equation*}
    E^{1}_{-2, q} = \cohom^{[q-2]}(X_2, Y_2) =
\begin{cases}
\CC^{l-1} & \text{if } q \text{ is odd}\\ 0 & \text{ otherwise}.
\end{cases}
\end{equation*}
The differential $d^1 : E^{1}_{-p, q} \to E^{1}_{-p-1, q}$ turns
out to be surjective for $p = 1$ and $q$ odd. For some obvious
geometric reasons, the spectral sequence $E^r_{p, q}$ at $r = 2$
then collapses and gives
\begin{equation}
    \Hp_{0}(A_1) \simeq E^{2}_{-1, 1} \simeq \CC^{l+1}
\end{equation}
and $\Hp_{1}(A_1) = 0$. We shall look again at the algebra $A_2$,
from a different point of view, in Example \ref{ex3}.
\end{example}

Let us consider now the following related example.

\begin{example}\label{ex2}\
Let $A_2 = M_2(\CC[x]) \oplus \CC^l$. Then the standard filtration
of $A_2$ is $\mfI_0 = A_2$, $\mfI_1 = M_2(\CC[x])$, and $\mfI_2 =
0$. The center of $A_2/\mfI_1$ is $Z_1 \simeq \CC^{l}$. We have
that $X_1$ consists of $l$ points and $Y_1 = 0$. The center of
$A_2/\mfI_2$ is $Z_2 \simeq \CC[x] \oplus \CC^l$. Then $X_2 = \CC
\cup X_1$ and $Y_2 = X_1$, where $X_1 \cap \CC = \emptyset$.

The spectral sequence of Theorem \ref{Theorem.spectral.sq} then
becomes $E^{1}_{-p, q}=0$, unless $p = 1$ or $p = 2$, in which case we
get
\begin{equation*}
    E^{1}_{-1, q} = \cohom^{[q-1]}(X_1, Y_1) =
\begin{cases}
\CC^{l} & \text{if } q \text{ is odd}\\ 0 & \text{ otherwise}
\end{cases}
\end{equation*}
and
\begin{equation*}
    E^{1}_{-2, q} = \cohom^{[q-2]}(X_2, Y_2) =
\begin{cases}
0 & \text{if } q \text{ is odd}\\ \CC & \text{ otherwise}.
\end{cases}
\end{equation*}
The spectral sequence collapses at the $E^1$ term for geometric
reasons and hence
\begin{equation}
    \Hp_{0}(A_2) \simeq E^{1}_{-1, 1} \oplus E^{1}_{-2, 2}
    \simeq \CC^{l+1}
\end{equation}
and $\Hp_{1}(A_1) = 0$.
\end{example}

The algebras $A_1$ and $A_2$ in the above examples turned out to have
the same periodic cyclic homology groups. These algebras are simple,
but representative, of the finite type algebras arising in the
representation theory of reductive $p$-adic groups. Clearly, the
periodic cyclic homology groups of these algebras provide important
information on the structure of these algebras, but fails to
distinguish them. At a heuristical level, distinguishing between $A_1$
and $A_2$ is the same problem as distinguishing between square
integrable representations and supercuspidal representations. This
issue arises because both these types of representations provide
similar homology classes in $\HP_0$ (through the Chern characters of
the idempotents defining them). It is then an important question to
distinguish between these homology classes.

Let us complete our discussion with some related results on the
Hochschild homology of finite type algebras. We begin with the
Hochschild homology of certain cross product algebras.

Let $\Gamma$ be a finite group acting on a smooth complex
algebraic variety $X$. For any $\gamma \in \Gamma$, let us denote
by $X^\gamma \subset X$ the points of $X$ fixed by $\gamma$. Let
\begin{equation*}
    \Gamma_\gamma := \{g\in \Gamma, g \gamma = \gamma g\}
\end{equation*}
denote the centralizer of $\gamma$ in $\Gamma$. Let $C_\gamma$ be the
(finite, cyclic) subgroup generated by $\gamma$. There exists a
natural $\Gamma$-invariant map
\begin{equation}
    \hat X := \bigcup_{\gamma \in \Gamma}
    \{\gamma\} \times
    X^\gamma \times (\Gamma_\gamma/C_\gamma)  \to X,
\end{equation} given
simply by the projection onto the second component. This gives then rise
to a $\Gamma$-equivariant morphism $\maO(X) \to \maO(\hat X)$. Choose
a representative $x \in \Gamma$ from each conjugacy class $\langle x
\rangle$ of $\Gamma$ and denote by $m_x$ the number of elements in the
conjugacy class of $m_x$. Denote by $C_\gamma^*$ the dual of
$C_\gamma$, that is the set of multiplicative maps $\pi : C_\gamma \to
\CC^*$.  Recall that $\<\Gamma\>$ denotes the set of conjugacy classes
of $\Gamma$. We are finally ready to define the morphisms
\begin{equation*}
    \psi_{\gamma, \pi} : \maO(X) \rtimes \Gamma \to
    M_{m_x}(\maO(X^\gamma)),
\end{equation*}
where $\gamma$ runs through a system of
representatives of the conjugacy classes of $\Gamma$ and $\pi \in
C_\gamma^*$
\begin{multline}\label{def.psi.g}
    \psi = \bigoplus_{\langle \gamma \rangle \in \<\Gamma \>, \pi
    \in C_\gamma^*} \psi_{\gamma, \pi}: \maO(X) \rtimes \Gamma \to
    \maO(\hat X) \rtimes \Gamma \\ \simeq \bigoplus_{\langle
    \gamma \rangle \in \<\Gamma \>} M_{m_x}(\maO(X^\gamma) \otimes
    \CC[C_\gamma])  \simeq \bigoplus_{\langle \gamma \rangle \in
    \<\Gamma \>, \pi \in C_\gamma^*} M_{m_x}(\maO(X^\gamma)).
\end{multline}
The equations above give rise to a map
\begin{equation}\label{eq.def.phi}
    \phi = \bigoplus_{\langle \gamma \rangle \in \langle \Gamma
    \rangle} \phi_\gamma : \HH_q(\maO(X) \rtimes \Gamma) \to
    \bigoplus_{\langle \gamma \rangle \in \langle \Gamma \rangle}
    \cohom^q(X^\gamma),
\end{equation}
defined using $\HH_q(M_{m}(\maO(X^\gamma))) \simeq
\HH_q(\maO(X^\gamma)) \simeq \cohom^q(X^\gamma)$ and
\begin{equation}\label{eq.def.phi2}
    \phi_\gamma = \sum_{\pi \in C_\gamma^*}
    \frac{\overline{\pi(\gamma)}}{\#C_\gamma}
     \HH_{q}(\psi_{\gamma, \pi}) :
    \HH_q(\maO(X) \rtimes \Gamma) \to \Omega^q(X^\gamma).
\end{equation}
(This map was defined in a joint work in progress with J. Brodzki
\cite{BrodzkiNistor}.)

Then we have the following lemma, which is a particular case of
the Theorem \ref{thm.cross} below.

\begin{lemma}\label{lemma.FG}\
Assume that $X = \CC^n$ and that $\Gamma$ acts linearly on
$\CC^n$. Then the map $\phi$ of Equation \ref{eq.def.phi} defines
an isomorphism
\begin{equation*}
    \phi : \HH_q(\maO(X) \rtimes \Gamma) \to
    \bigoplus_{\langle \gamma \rangle \in \langle \Gamma \rangle}
    \Omega^q(X^\gamma)^{\Gamma_\gamma}.
\end{equation*}
\end{lemma}

\begin{proof}\
Let $A$ be a complex algebra acted upon by automorphisms by a group
$\Gamma$. Let us recall \cite{FeiginTsygan1, Nistor8} that the groups
$\HH_q(A \rtimes \Gamma)$ decompose naturally as a direct sum
\begin{equation}
    \HH_q(A \rtimes \Gamma) \simeq \oplus_{\langle \gamma \rangle \in \langle
    \Gamma \rangle} \HH_q(A \rtimes \Gamma)_{\gamma}.
\end{equation}
The components $\HH_q(A \rtimes \Gamma)_{\gamma}$ are then
identified as follows. Let $A_\gamma$ be the $A$--$A$ bimodule
with action $a \cdot b \cdot c = ab\gamma(c)$. Let $\HH_q(A, M)$
denote the Hochschild homology groups of an $A$--bimodule $M$
\cite{MacLane, Moerdijk, Weibel}. Then
\begin{equation}\label{eq.1}
    \HH_q(A \rtimes \Gamma)_{\gamma}
    \simeq \HH_q(A, A_\gamma)^{\Gamma_\gamma}.
\end{equation}
This follows from example from \cite{Nistor8}[Lemma 3.3] (take $G
= \Gamma$ in that lemma). This will also be discussed in
\cite{BrodzkiNistor, Shantanu}.

Let now $A = \maO(X)$, with $X = \CC^n$ and $\Gamma$ acting
linearly on $X$. The same method as the one used in the proof of
\cite{BrylinskiNistor}[Lemma 5.2] shows that
\begin{equation}\label{eq.2}
    \HH_q(A, A_\gamma) \simeq \HH_q(\maO(X^{\gamma}),
    \maO(X^{\gamma})) = \HH_q(\maO(X^{\gamma})) =
    \Omega^q(X^{\gamma}),
\end{equation}
the isomorphism being given by the restriction morphism $\maO(X)
\to \maO(X^\gamma)$.

A direct calculation based on the formula for the map $J$ in
\cite{Nistor8}[Lemma 3.3] shows that the composition of the
morphisms of Equations \eqref{eq.1} and \eqref{eq.2} is the map
$\phi_\gamma$ of Equation \eqref{eq.def.phi2}. This completes the
proof.
\end{proof}

We now extend the above result to $X$ an arbitrary smooth, complex
algebraic variety. To do this, we shall need two results on
Hochschild homology. We begin with the following result of
Brylinski from \cite{Brylinski1}.

\begin{proposition}[Brylinski]\label{Prop.localization} \
Let $S$ be a multiplicative subset of the center $Z$ of the
algebra $A$. Then $\Hd_*(S^{-1}A) \simeq S^{-1}\Hd_*(A)$.
\end{proposition}

A related result, from \cite{KNS}, studies the completion of
Hochschild homology. In the following result, we shall use {\em
topological Hochschild homology}, whose definition is similar to
that of the usual Hochschild homology, except that one completes
with respect to the powers of an ideal (below, we shall use this
for the ideals $I$ and $IA$).

\begin{theorem}\label{theorem.completion} \
Suppose that $A$ is a unital finite type $\kk$--algebra and $I
\subset \kk$ is an ideal. Then the natural map $\Hd_*(A) \to
\tHd_*(\hat A)$ and the $\kk$--module structure on $\Hd_*(A)$
define an isomorphism
\begin{equation*}
    \Hd_*(A) \otimes_\kk \hat \kk \simeq \tHd_*(\hat A)
\end{equation*}
of $\hat \kk$--modules.
\end{theorem}

We are ready now to prove the following theorem.

\begin{theorem}\label{thm.cross}\ Assume that $X$ is a smooth,
complex algebraic variety and that $\Gamma$ acts on $X$ by algebraic
automorphisms. Then the map $\phi$ of Equation \ref{eq.def.phi}
defines an isomorphism
\begin{equation*}
    \phi : \HH_q(\maO(X) \rtimes \Gamma) \to
    \bigoplus_{\langle \gamma \rangle \in \langle \Gamma \rangle}
    \Omega^q(X^\gamma)^{\Gamma_\gamma}.
\end{equation*}
\end{theorem}

\begin{proof}\ Let
\begin{equation*}
    Z := \maO(X)^{\Gamma} = \maO(X/\Gamma).
\end{equation*}
Then the morphism $\psi$ of Equation \eqref{def.psi.g} is
$Z$-linear. It follows that the map $\phi$ of Equation
\eqref{eq.def.phi} is also $Z$-linear. It is enough hence to prove
that the localization of this map to any maximal ideal of $Z$ is
an isomorphism. It is also enough to prove that the completion of
any of these localizations with respect to that maximal ideal is
an isomorphism. Since $X$ is smooth (and hence the completion of
the local ring of $X$ at any point is a power series ring) this
reduces to the case of $X = \CC^n$ acted upon linearly by
$\Gamma$. The result hence follows from Lemma \ref{lemma.FG}.
\end{proof}

Let us include the following corollary of the above proof.

\begin{corollary}\ The map $\phi_\gamma : \HH_{q}(\maO(X) \rtimes \Gamma)
\to \Omega^q(X^\gamma)$ of \eqref{eq.def.phi2}
is such that $\phi_\gamma = 0$ on $\HH_{q}(\maO(X) \rtimes \Gamma)_{g}$
if $g$ and $\gamma$ are not in the same conjugacy class and induces
an isomorphism
\begin{equation*}
    \phi_\gamma : \HH_{q}(\maO(X) \rtimes \Gamma)_{\gamma}
    \to \Omega^q(X^\gamma)^{\Gamma_\gamma}.
\end{equation*}
\end{corollary}

Let us apply these results to the algebra $A_1$ of Example \ref{ex1}.

\begin{example}\label{ex.1prime}\
Let $\Gamma = \ZZ/2\ZZ$ act by $z \to -z$ on $\CC$.  Chose $l = 1$ and
$a_1 = 0$ in Example \ref{ex1}. Then $\maO(\CC) \rtimes \Gamma \simeq
A_1$, and hence
\begin{equation*}
    \HH_q(A_1) \simeq
    \begin{cases}
    \maO(\CC)^{\Gamma} \oplus \CC & \text{ if } q = 0 \\
    \maO(\CC)^{\Gamma} & \text{ if } q = 1 \\
    0 & \text{ otherwise.}
    \end{cases}
\end{equation*}
In particular, $\HH_q(A_1)$ vanishes for $q$ large. Also, note
that
\begin{equation*}
\maO(\CC)^{\Gamma} = \CC[x]^{\Gamma} = \CC[x^2] \simeq \CC[x].
\end{equation*}
\end{example}

Let us see now an example of a finite type algebra for which
Hochschild homology does not vanish in all degrees.

\begin{example}\label{ex.nb}\
Let $A_3 \subset M_2(\maO(\CC)) = M_2(\CC[x])$ be the subalgebra of
those matrices
\begin{equation*}
P = \left [
\begin{array}{cc}
P_{11} & P_{12} \\
P_{21} & P_{22}
\end{array}
\right ] =: P_{11}e_{11} +  P_{12} e_{12}
P_{21}e_{21} +  P_{22},
\end{equation*}
with the property that $P_{21}(0) = P_{21}'(0) = 0$.  For a
suitable choice of $v_1$, this is a subalgebra of the algebra
$A_1$ considered in Example \ref{ex.1prime}. Let $V_1 :=
A_3e_{11}$ and $V_2 := A_3e_{22}$. Then $M := V_2/V_1 \simeq
\CC[x]/(x^2)$. The modules $V_1$ and $V_2^{\tau} := e_{22}A_3$ can
be used to produced a projective resolution of $A_3$ with free
$A_3$--$A_3$ bimodules that gives
\begin{equation}
    \HH_q(A_3) \simeq
    \begin{cases}
    \CC[x] \oplus M & \text{ if } q = 0 \\
    \CC[x] \oplus \CC & \text{ if } q = 1 \\
    \CC & \text{ otherwise.}
    \end{cases}
\end{equation}
\end{example}

\section{Spectrum preserving morphisms\label{sec.BC}\label{sec3}}

We shall now give more evidence for the close relationship between the
topology of $\Prim(A)$ and $\Hp_*(A)$ by studying a class of morphisms
implicitly appearing in Lusztig's work on the representation of
Iwahori-Hecke algebras, see \cite{LC1, LC2, LC3, LC4}.

Let $L$ and $J$ be two finite type $\kk$-algebras. If $\phi : L
\to J$ is a $\kk$-linear morphism, we define
\begin{equation}\label{eq.sp.rel}
    {\mathcal R}_\phi := \{ (\mfk P',\mfk P) \subset \Prim(J)
    \times \Prim(L),\, \phi^{-1}(\mfk P') \subset \mfk P \}.
\end{equation}

We now introduce the class of morphisms we are interested in.

\begin{definition}\label{Def.mult.one}\
Let $\phi : L \to J$ be a $\kk$-linear morphism of unital, finite
type $\kk$-algebras. We say that $\phi$ is a {\em spectrum
preserving morphism} if, and only if, the set $\mathcal R_\phi$
defined in Equation \eqref{eq.sp.rel} is the graph of a bijective
function
\begin{equation*}
    \phi^* : \Prim(J) \to \Prim(L).
\end{equation*}
\end{definition}

More concretely, we see that $\phi : L \to J$ is spectrum
preserving if, and only if, the following two conditions are
satisfied:

\begin{enumerate}
\item\ {\em For any primitive ideal $\mfk P$ of $J$, the ideal
$\phi^{-1}(\mfk P)$ is contained in a unique primitive ideal of
$L$, namely $\phi^*(\mfk P)$}, and
\item\ {\em The resulting map $\phi^*: \Prim(J) \to \Prim(L)$ is a
bijection}.
\end{enumerate}

We have the following result combining two theorems from
\cite{BaumNistor}.

\begin{theorem}[Baum-Nistor]\
Let $L$ and $J$ be finite type $\kk$--algebras and $\phi : L \to
J$ be a $\kk$--linear spectrum preserving morphism. Then the
induced map $\phi^* : \Prim(J) \to \Prim(L)$ between primitive
ideal spectra is a homeomorphism and the induced map $\phi_* :
\Hp_*(L) \to \Hp_*(J)$ between periodic cyclic homology groups is
an isomorphism.
\end{theorem}

We also obtain an isomorphism on periodic cyclic homology for a
slightly more general class of algebra morphisms.

\begin{definition}\label{def.w.sp.p}\
A morphism $\phi : L \to J$ of finite type algebras is called {\em
weakly spectrum preserving} if, and only if, there exist
increasing filtrations
\begin{equation*}
\begin{gathered}
    (0) = L_0 \subset L_1 \subset \ldots \subset L_n = L \quad
    \text{ and } \qquad (0) = J_0 \subset J_1 \subset \ldots
    \subset J_n = J
\end{gathered}
\end{equation*}
of two-sided ideals such that $\phi(L_k) \subset J_k$ and the
induced morphisms $L_k/L_{k-1} \to J_k / J_{k-1}$ are spectrum
preserving.
\end{definition}

Combining the above theorem with the excision property, we obtain
the following result from \cite{BaumNistor}.

\begin{theorem}[Baum-Nistor]\
Let $L$ and $J$ be finite type $\kk$--algebras and $\phi : L \to
J$ be a $\kk$--linear weakly spectrum preserving morphism. Then
the induced map $\phi_* : \Hp_*(L) \to \Hp_*(J)$ between periodic
cyclic homology groups is an isomorphism.
\end{theorem}

The main application of this theorem is the determination of the
periodic cyclic homology of Iwahori-Hecke algebras in
\cite{BaumNistor0, BaumNistor}. Let $H_q$ be an Iwahori-Hecke
algebra and $J$ the corresponding asymptotic Hecke algebra
associated to an extended affine Weyl group $\widehat{W}$
\cite{LC3} (their definition is recalled in \cite{BaumNistor},
where more details and more complete references are given). Then
there exists a morphism $\phi: H_q \to J$ of $\kk$--finite type
algebras, $\kk = Z(H_q)$ that is weakly spectrum-preserving
provided that $q$ is not a root of unity or $q = 1$. Therefore
\begin{equation}
    \phi_* : \Hp_*(H_q) \to \Hp_*(J)
\end{equation}
is an isomorphism. The algebra $H_1$ is a group algebra, and hence
its periodic cyclic homology can be calculated directly.

The above theorem also helps us clarify the Examples \ref{ex1} and
\ref{ex2}.

\begin{example}\label{ex3}\
Let us assume that $v_j = v \in \CC^2$ in Example \ref{ex1}. Let
$e \in M_2(\CC)$ be the projection onto the vector $e$. The
morphism $\phi: \CC[x] \ni P \to Pe \in A_1 \subset M_2(\CC[x])$
is not weakly spectrum preserving. However, $\phi : \CC[x] \to
\mfI_1$ is a spectrum preserving morphism of $\CC[x]$--algebras.
Combining with the inclusion $\mfI_1 \subset M_2(\CC[x])$, we see
that $\phi_* : \Hp_*(\CC[x]) \to \Hp_*(\mfI_1)$ is an isomorphism
and $\Hp_*(\mfI_1)$ is a direct summand of $\Hp_*(A_1)$. The
excision theorem then gives
\begin{equation*}
    \Hp_q(A_1) \simeq \Hp_q(\mfI_1) \oplus \Hp_q(A/\mfI_1)
    \simeq \Hp_q(\CC[x]) \oplus \Hp_q(\CC^l) =
    \begin{cases}
    \CC^{l+1} & \text{if } q \text{ is even}\\
    0 & \text{ otherwise.}
    \end{cases}
\end{equation*}
\end{example}

The case of the algebra $A_2$ of example is even simpler.

\begin{example}\label{ex4}\
The inclusion $Z(A_2) \to A_2$ is a spectrum preserving morphism
of $Z(A_2)$--algebras. Consequently,
\begin{equation*}
    \Hp_q(A_2) \simeq \Hp_q(\CC[x]) \oplus \Hp_q(\CC^l) =
    \begin{cases}
    \CC^{l+1} & \text{if } q \text{ is even}\\
    0 & \text{ otherwise.}
    \end{cases}
\end{equation*}
\end{example}

Let us notices that by considering the action of the natural
morphisms
\begin{equation*}
\CC[x] \to Z(A_1) \subset A_1,\;  A_1 \to M_2(\CC[x]),\; \CC[x]
\to Z(A_2) \subset A_2,\; or A_2 \to M_2(\CC[x])
\end{equation*}
on periodic cyclic homology, we will still not be able to distinguish
between $A_1$ and $A_2$. However, the natural products $\Hp_i(\CC[x])
\otimes \Hp_j(A_k) \to \Hp_*(A_k)$ (see \cite{Kassel, Kassel1}) {\em
will} distinguish between these algebras.

\section{The periodic cyclic homology of $\CIc(G)$\label{sec4}}

Having discussed the relation between $\Hp_*(\CIc(G))$ and the
admissible spectrum $\hat{G} = \Prim(\CIc(G))$, let us recall the
explicit calculation of $\Hp_*(\CIc(G))$ from \cite{Nistor29}. The
calculation of $\Hp_*(\CIc(G))$ in \cite{Nistor29} follows right away
from the calculation of the Hochschild homology groups of
$\CIc(G)$. The calculations of presented in this section complement
the results on the cyclic homology of $p$--adic groups in
\cite{Nistor14, Schneider}.

To state the main result of \cite{Nistor29} on the Hochschild
homology of the algebra $\CIc(G)$, we need to introduce first the
concepts of a ``standard subgroup'' and of a ``relatively regular
element'' of a standard subgroup. For any group $G$ and any subset
$A \subset G$, we shall denote
\begin{equation*}
    C_G(A) := \{ g \in G, ga = ag, \ \forall a \in A \}\,,\quad
    N_G(A) := \{ g \in G, gA = Ag\}\,,
\end{equation*}
$W_G(A) := N_G(A)/ C_G(A)$, and $Z(A) := A \cap C_G(A)$. This
latter notation will be used only when $A$ is a subgroup of $G$.
The subscript $G$ will be dropped from the notation whenever the
group $G$ is understood.  A commutative subgroup $S$ of $G$ is
called {\em standard} if $S$ is the group of semisimple elements
of the center of $C(s)$ for some semi-simple element $s \in G$. An
element $s \in S$ with this property will be called {\em regular
relative to $S$}, or {\em $S$-regular}. The set of $S$-regular
elements will be denoted by $S^{\rg}$.

We fix from now on a $p$-adic group $G$. (Recall that a $p$--adic
group $G= \GG(\FF)$ is the set of $\FF$--rational points of a
linear algebraic group $\GG$ defined over a non-archemedean,
non-discrete, locally compact field $\FF$ of characteristic zero.)
Our results will be stated in terms of standard subgroups of $G$.
We shall denote by $H_u$ the set of unipotent elements of a
subgroup $H$.  Sometimes, the set $C(S)_u$ is also denoted by
$\nil S$, in order to avoid having to many paranthesis in our
formulae. Let $\Delta_{C(S)}$ denote the modular function of the
group $C(S)$ and let
\begin{equation*}
    \CIc(\nil S)_{\delta}:=\CIc(C(S)_u) \otimes \Delta_{C(S)},
\end{equation*}
be $\CIc(\nil S)$ as a vector space, but with the product
$C(S)$-module structure, that is
\begin{equation*}
    \gamma(f)(u) = \Delta_{C(S)}(\gamma) f(\gamma^{-1}u\gamma),
\end{equation*}
for all $\gamma \in C(S)$, $f \in \CIc(\nil S)_\delta$ and $u \in \nil
S$.

The groups $\Hd_*(\CIc(G))$ are determined in terms of the
following data:
\begin{enumerate}
\item\ the set $\Sigma$ of conjugacy classes of standard subgroups $S$
of $G$;
\item\ the subsets $S^{\rg} \subset S$ of $S$-regular elements;
\item\ the actions of the Weyl groups $W(S)$ on $\CIc(S)$;\ and
\item\ the continuous cohomology of the $C(S)$--modules $\CIc(\nil
S)_{\delta}$.
\end{enumerate}

Combining this proposition with Corollary \ref{thm.DS}, we obtain
the main result of this section. Also, recall that $\nil S$ is the
set of unipotent elements commuting with the standard subgroup
$S$, and that the action of $C(S)$ on $\CIc(\nil S)$ is twisted by
the modular function of $C(S)$, yielding the module $\CIc(\nil
S)_{\delta}= \CIc(\nil S) \otimes \Delta_{C(S)}$.

\begin{theorem} \label{Theorem.Gen}\
Let $G$ be a $p$--adic group. Let $\Sigma$ be a set of
representative of conjugacy classes of standard subgroups of $S
\subset G$ and $W(S) = N(S)/C(S)$, then we have an isomorphism
\begin{equation*}
    \Hd_q(\CIc(G)) \simeq \bigoplus\limits_{S \in \Sigma}
    \CIc(S^{\rg})^{W(S)} \otimes \cohom_q(C(S), \CIc(\nil S)_{\delta}).
\end{equation*}
\end{theorem}

The isomorphism of this theorem was obtained by identifying the
$E^\infty$-term of a spectral sequence convergent to
$\Hd_q(\CIc(G))$, and hence it is not natural. This isomorphism
can be made natural by using a generalization of the Shalika germs
\cite{Nistor29} .

The periodic cyclic homology groups of $\CIc(G)$ are then
determined as follows. Recall that our convention is that
$\Hd_{[q]} := \oplus_{k \in \ZZ} \Hd_{q + 2k}$. An element $\gamma
\in G$ is called {\em compact} if the closure of the set
$\{\gamma^n\}$ is compact. The set $G_{comp}$ of compact elements
of $G$ is open and closed in $G$, if we endow $G$ with the locally
compact, Hausdorff topology obtained from an embedding $G \subset
\FF^N$. We clearly have $\gamma G_{comp}\gamma^{-1} = G_{comp}$,
that is, $G_{comp}$ is $G$-invariant for the action of $G$ on
itself by conjugation. Also, we shall denote by
$\Hd_{[q]}(\CIc(G))_{comp}$ the localization of the homology group
$\Hd_{[q]}(\CIc(G))$ to the set of compact elements of $G$ (see
\cite{Blanc-Brylinski} or \cite{Nistor8}). This localization is
defined as follows. Let $R^{\infty}(G)$ be the ring of locally
constant $Ad_G$-invariant functions on $G$ with the pointwise
product. If $\omega = f_0 \otimes f_1 \otimes \ldots f_n \in
\CIc(G)^{\otimes (n + 1)} = \CIc(G^{n+1})$ and $z \in
R^{\infty}(G)$, then we define
\begin{equation*}
    [z \omega](g_0, g_1, \ldots, g_n) = z(g_0g_1 \ldots g_n)
    \omega (g_0, g_1, \ldots, g_n) \in \CIc(G^{n+1}).
\end{equation*}
Let $\chi$ be the characteristic function of $G_{comp}$ (so $\chi
= 1$ on $G_{comp}$ and $\chi = 0$ otherwise). Then
$\Hd_{[q]}(\CIc(G))_{comp} = \chi \Hd_{[q]}(\CIc(G))$.

\begin{theorem}\ We have
\begin{equation}\label{eq.HN}
    \Hp_q(\CIc(G)) \simeq \HH_{[q]}(\CIc(G))_{comp}.
\end{equation}
Let $S_{comp}$ be the set of compact elements of a standard
subgroup $S$, then
\begin{equation} \label{eq.result2}
    \Hp_q (\CIc(G)) \simeq \bigoplus \limits_{S \in \Sigma}
    \CIc(S_{comp}^{\rg})^{W(S)} \otimes \cohom_{[q]} (C(S),
    \CIc(\nil S)_\delta).
\end{equation}
\end{theorem}

The Equation \eqref{eq.HN} follows also from the results in
\cite{Nistor14}.

It is conceivable that a next step would be to study the
``discrete parts'' of the groups $\HH_q(\CIc(G))$, following the
philosophy of \cite{BDK, Kazhdan}. This can be defined as follows.
In \cite{Nistor29}, we have defined ``induction morphisms''
\begin{equation}
    \phi_M^G : \HH_q(\CIc(G)) \to \HH_q(\CIc(M))
\end{equation}
for every Levi subgroup $M \subset G$. We define
$\HH_q(\CIc(G))_{0}$ to be the intersection of all kernels of
$\phi_M^G$, for $M$ a proper Levi subgroup of $G$ and call this
the discrete part of $\HH_q(\CIc(G))_{0}$.

Assuming that one has established, by induction, a procedure to
construct the cuspidal representations of all $p$--adic groups of
lower split rank, then one can study the action of the center of
$\CIc(G)$ corresponding to the cuspidal associated to proper Levi
subgroups on the discrete part of $\HH_q(\CIc(G))$, which would
hopefully allow us to distinguish between the cuspidal part of
$\HH_q(\CIc(G))_{0}$ and its part coming from square integrable
representations that are not super-cuspidal.

\section{Hochschild homology \label{sec5}}

We include now three short sections that recall some of the
definitions used above. Nothing in this and the next section is
new, and the reader interested in more details as well as precise
references should consult the following standard references
\cite{Brodzki, BBN, ConnesNCG, ConnesBook, Karoubi, LodayQuillen,
Loday, Tsygan}.

We begin by recalling the definitions of Hochschild homology
groups of a complex algebra $A$, not necessarily with unit. We
define $b$ and $b'$ define two linear maps
\begin{equation}
    b,\, b'\, :\, A^{\otimes n + 1} \to A^{\otimes n}.
\end{equation}
where $A^{\otimes n} := A \otimes A \otimes \ldots A$\ ($n$ times)
by the formulas
\begin{equation}\label{eq.def.bb'}
\begin{gathered}
    b'(\Tt) = \sum_{i=0}^{n-1} (-1)^ia_0\otimes\ldots\otimes a_i
    a_{i+1}\otimes\ldots\otimes a_n,\\ b(\Tt)=b'(\Tt) +(-1)^n
    a_na_0\otimes\ldots\otimes a_{n-1},
\end{gathered}
\end{equation}
where $a_0, a_1, \ldots, a_n \in A$. Let
\begin{equation}
    \mathcal{H}_n(A) = \mathcal{H}_n'(A) := A \otimes A \otimes
    \ldots A \; (n+1 \text{ times }).
\end{equation}
Also, let $\maH(A) := (\maH_n(A), b)$ and $\maH'(A) :=
(\maH/_n(A), b')$. The homology groups of $\maH(A)$ are, by
definition, the {\em Hochschild homology} groups of $A$ and are
non-zero. The $n$th Hochschild homology group of $A$ is denoted
$\HH_n(A)$. By dualizing, we obtain the Hochschild cohomology
groups $\HH^n(A)$.

If $A$ has a unit, then the complex $\maH'(A)$ is acyclic (\ie it
has vanishing homology groups) because $b's + s b' = 1$, where
\begin{equation}\label{eq.def.s}
    s(\Tt)=1\otimes \Tt.
\end{equation}
Therefore, if $A$ has a unit, the complex $\mathcal{H}'(A)$ is a
resolution of $A$ by free $A$-bimodules.

Recall that {\em a trace} on $A$ is a linear map $\tau : A \to
\CC$ such that $\tau(a_0 a_1 - a_1 a_0) = 0$ for all $a_0, a_1 \in
A$. The space of all traces on $A$ is then isomorphic to
$\Hd^0(A)$.

An algebra $A$ such that $\mathcal{H}'(A)$ is acyclic is called
{\em $H$-unital}, following Wodzicki \cite{Wodzicki}.

Clearly the groups $\Hd_n(A)$ are {\em covariant} functors in $A$,
in the sense that any algebra morphism $\phi: A \to B$ induces a
morphism
\begin{equation*}
    \phi_* = \Hd_n(\phi) : \Hd_n(A) \to \Hd_n(B)
\end{equation*}
for any integer $n\geq 0$. Similarly, we also obtain a morphism
$\phi^* = \Hd^n(\phi): \Hd^n(B) \to \Hd^n(A)$. In other words,
Hochschild cohomology is a {\em contravariant} functor.  It is
interesting to note that if $Z$ is the center of $A$, then $\Hd_n(A)$
is also a $Z$-module, where, at the level of complexes the action is
given by
\begin{equation}\label{eq.module}
    z(\Tt) = za_0\otimes a_1 \otimes \ldots \otimes a_n.
\end{equation}
for all $z\in Z$. As $z$ is in the center of $A$, this action will
commute with the Hochschild differential $b$.

\section{Cyclic homology\label{sec6}}

Let  $A$ be a unital algebra. We  shall denote by $t$ the (signed)
generator of cyclic permutations:
\begin{equation}
    t(\Tt) = (-1)^n a_n\otimes
    a_0 \otimes \ldots \otimes a_{n-1}
\end{equation}
Using the operator $t$ and the contracting homotopy $s$ of the complex
$\mathcal{H}'(A)$, Equation \eqref{eq.def.s}, we construct a new
differential $B := (1 - t)B_0$, of degree $+1$, where
\begin{equation}
    B_0(\Tt) = s\sum_{k=0}^{n} t^k(\Tt).
    \end{equation}
It is easy to check that $B^2 = 0$ and that  $[b,B]_+ := bB + Bb =  0$.

The differentials $b$ and $B$ give rise to the following complex
\begin{figure}[ht] \setlength{\unitlength}{1ex}
\begin{picture}(60,28)(0,12)

\put(15,38){\vector(0,-1){6.2}} \put(30,38){\vector(0,-1){6.2}}
\put(45,38){\vector(0,-1){6.2}} \put(12.5,35){\makebox{$b$}}
\put(30.6,35){\makebox{$b$}} \put(46,35){\makebox{$b$}}
\put(12.5,29){\makebox(5,1){$A^{\otimes 3} $}}
\put(25.2,29.5){\vector(-1,0){6}}
\put(27.5,29){\makebox(5,1){$A^{\otimes 2}$}}
\put(40.2,29.5){\vector(-1,0){6}}
\put(42.5,29){\makebox(5,1){$A$}}
\put(22.5,30.5){\makebox{$B$}} \put(37.2,30.5){\makebox{$B$}}
\put(15,27.3){\vector(0,-1){5.9}} \put(30,27.3){\vector(0,-1){5.9}}
\put(12.5,24.3){\makebox{$b$}} \put(30.6,24.3){\makebox{$b$}}
\put(12.7,19){\makebox(5,1){$A^{\otimes 2} $}}
\put(26.5,19.5){\vector(-1,0){8}}
\put(27.7,19){\makebox(5,1){$A$}}

\put(22.5,20.5){\makebox{$B$}}

\put(15,17.5){\vector(0,-1){6}}  \put(12.5,14.5){\makebox{$b$}}

\put(12.7,9){\makebox(5,1){$A$}} \end{picture}
\vspace{.2in} \caption{The cyclic bicomplex of the algebra $A$.}
\end{figure}

We notice that columns the above complex are copies of the
Hochschild complex $\mathcal{H}(A)$. The cyclic complex $\mathcal
C(A)$ is by definition the total complex of this double complex.
Thus the space of cyclic $n$-chains is defined by
\begin{equation}
    {\mathcal C}(A)_n = \bigoplus_{k\geq 0}
    {\mathcal H}_{n - 2k}(A),
\end{equation}
we see that $({\mathcal C}(A), b+B),$ is a complex, called {\em
the cyclic complex} of $A$, whose homology is by definition the
{\em cyclic homology} of $A$, denoted $\Hc_q(A)$, $q \ge 0$.

There is a canonical operator $S: \mathcal{C}_{n}(A) \rightarrow
\mathcal{C}_{n-2}(A)$, called the \emph{Connes periodicity
operator}, which shifts the cyclic complex left and down,
explicitly defined by
\begin{equation}
    S(\omega_n, \omega_{n-2}, \omega_{n-4}, \dots ) \mapsto
(\omega_{n-2}, \omega_{n-4}, \dots ),
\end{equation}
where $\omega_k\in \mathcal{H}_k(A)$, for all $k$.  This operator
induces the short exact sequence of complexes
\begin{equation*}
    0 \rightarrow {\mathcal H}(A) \xrightarrow{I} {\mathcal C}(A)
    \xrightarrow{S} {\mathcal C}(A)[2]\rightarrow 0,
\end{equation*}
where the map $I$ is the inclusion of the Hochschild complex as the
first column of the cyclic complex. The snake Lemma in homology
\cite{MacLane, Moerdijk, Weibel} gives the following long exact
sequence, called the {\em SBI--exact sequence}
\begin{equation} \label{eq.SBI}
    \ldots\rightarrow \Hd_{n}(A)\stackrel{I}{\longrightarrow}
    \Hc_{n}(A)\stackrel{S}{\longrightarrow}
    \Hc_{n-2}(A)\stackrel{B}{\longrightarrow}
    \Hd_{n-1}(A)\stackrel{I}{\longrightarrow}\ldots\, ,
\end{equation}
where $B$ is the differential defined above, see \cite{ConnesNCG,
LodayQuillen} for more details.  The {\em periodic cyclic complex}
of an algebra $A$ is the complex
\begin{equation*}
    \maC^{\per}(A) := \displaystyle{\lim_{\leftarrow}} \, \maC(A),
\end{equation*}
the inverse limit being taken with respect to the periodicity
morphism $S$. It is a $\ZZ/2\ZZ$-graded complex, whose chains are
(possibly infinite) sequences of Hochschild chains with degrees of
the same parity. The homology groups of the periodic cyclic
complex $\maC^{\per}(A)$ are, by definition, the {\em periodic
cyclic homology groups} of $A$. A simple consequence of the
SBI--exact sequence is that if $\phi : A \to B$ is a morphism of
algebras that induces an isomorphism on Hochschild homology, then
$\phi$ induces an isomorphism on cyclic and periodic cyclic
homology groups as well. Here is an application of this simple
principle.

Here is an application of this lemma. Consider
\begin{equation}
    Tr_* : \maH_q(M_N(A)) \to \maH_q(A), \quad q \in \ZZ_+,
\end{equation}
the map defined by $Tr_*(b_0 \otimes \ldots \otimes b_q) =
Tr(m_0m_1 \ldots m_q) a_0 \otimes \ldots a_q$, if $b_k = m_k
\otimes a_k \in M_N(\CC) \otimes A = M_N(A)$. Also consider the
(unital) inclusion $\iota : A \to M_N(A)$ and $\iota_*$ be the
morphism induced on the Hochschild complexes.

\begin{proposition}\label{prop.Morita}\ The map $Tr_*$
commutes with $b$ and $B$. Both $\iota_*$ and $\Tr_*$ induce
isomorphisms on Hochschild, cyclic, and periodic cyclic homologies and
cohomologies such that $(\iota_*)^{-1} = N^{-1}Tr_*$.
\end{proposition}

The cyclic and periodic cyclic homology groups of a non-unital
algebra $A$ are defined as the cokernels of the maps $\Hc_q(\CC)
\to \Hc_q(A^+)$ and $\Hp_q(\CC) \to \Hp_q(A^+)$, where $A^+ = A
\oplus \CC$ is the algebra with adjoined unit.

More generally, cyclic homology groups can be defined for ``mixed
complexes,'' \cite{BBN, JonesKassel1, Kassel1, KNS}.

\section{The Chern character\label{sec7}}

We shall use these calculations to construct Chern characters. By
taking $X$ to be a point in Theorem \ref{thm.FT1}, we obtain that
$\Hc_{2q}(\CC) \simeq \CC$ and $\Hc_{2q + 1}(\CC) \simeq 0$. We
can take these isomorphisms to be compatible with the periodicity
operator $S$ and such that for $q = 0$ it reduces to
\begin{equation*}
    \Hc_{0}(\CC) = \Hd_{0}(\CC) = \CC/[\CC,\CC] = \CC.
\end{equation*}
We shall denote by $\eta_q \in \Hc_{2q}(\CC)$ the unique element
such that
\begin{equation}
    S^q \eta_q = 1 \in \CC = \Hc_{0}(\CC).
\end{equation}

For any projection $e \in M_N(A)$ is a projection we obtain a
(non-unital) morphism $\psi : \CC \to M_N(A)$ by $\lambda\mapsto
\lambda e$. Then {\em Connes-Karoubi Chern character of $e$ in
cyclic homology} \cite{ConnesNCG, ConnesBook, Karoubi, Loday} is
defined by
\begin{equation}
    Ch_q([e]) = Tr_* (\psi_*(\eta_q)) \in \Hc_{2q}(A).
\end{equation}
This map can be shown to depend only on the class of $e$ in $K$-theory
and to define a morphism
\begin{equation}
    Ch_q : K_0(A) \to \Hc_{2q}(A).
\end{equation}

One can define similarly the Chern character in periodic cyclic
homology and the Chern character on $K_1$ (algebraic $K$-theory).
For the Connes-Karoubi Chern character on $K_1$, we use instead $Y
= \CC^*$, whose algebra of regular functions is $\maO[Y] \simeq
\CC[z, z^{-1}]$, the algebra of Laurent polynomials in $z$ and
$z^{-1}$ (this algebra, in turn, is isomorphic to the group
algebra of $\ZZ$). Then $ \Hc_q(\maO[Y]) \simeq \CC$, for any $q
\ge 1$.  We are interested in the odd groups, which will be
generated by elements $v_k \in \Hc_{2k + 1}(\maO[Y])$, which can
be chosen to satisfy $v_1 = z^{-1} \otimes z$ and $S^kv_{2k + 1} =
v_1$.

Then, if $u \in M_N(A)$ is an invertible element, it defines a
morphism $\psi : \CC[\CC^*] \to M_N(A)$. The {\em Connes-Karoubi
Chern character of $u$ in cyclic homology} is thus defined by
\begin{equation}
    Ch_q([u]) = Tr_* (\psi_*(v_q)) \in \Hc_{2q + 1}(A).
\end{equation}
Again, this map can be shown to depend only on the class of $u$ in
$K$-theory and to define a morphism
\begin{equation}
    Ch_q : K_1(A) \to \Hc_{2q + 1}(A).
\end{equation}

Both the Chern character on $K_0$ and on $K_1$ are functorial, by
construction.



\end{document}